\documentclass[12pt,journal,onecolumn ]{IEEEtran}

\IEEEoverridecommandlockouts              

\usepackage{url}
\usepackage{amsmath,amssymb}
\usepackage{amsthm,thmtools}
\usepackage{amsfonts} 

\usepackage{enumerate}
\usepackage{color}                  

\usepackage{verbatim}

\usepackage{amssymb}
\usepackage{amsbsy}
\usepackage{amsmath,lipsum}
\usepackage{setspace}
\usepackage{url}

\usepackage{float}
 \usepackage{dsfont}
 \usepackage{cuted} \stripsep-3pt

\newtheorem{Problem}{Problem}

 \newcommand{\Int}{\operatorname{int}}
 
 \newcommand{\real}{\operatorname{Re}}
 
 \newcommand{\diag}{\operatorname{diag}}

 \newcommand{\adju}{\operatorname{adj}}
 
\newcommand{\st}{\, | \,}

\newcommand*\diff{\mathop{}\!\mathrm{d}}

\newcommand{\trace}{\operatorname{trace}}

\declaretheorem[name={Example},qed={\lower-0.3ex\hbox{$\square$}} ] {Example}

\declaretheorem[name={Definition}  ] {Definition}
\declaretheorem[name={Theorem}  ] {Theorem}
\declaretheorem[name={Lemma}  ] {Lemma}
\declaretheorem[name={Remark}  ] {Remark}
\declaretheorem[name={Corollary}  ] {Corollary}
\declaretheorem[name={Assumption}  ] {Assumption}
\declaretheorem[name={Proposition}  ] {Proposition}

\newcommand {\R}{\mathbb R}

\newcommand {\C}{\mathbb C}

\newcommand {\Q}{\mathbb Q}

\renewcommand {\P}{\mathbb P}

\newcommand {\EP}{\mathbb E \mathbb P}

\newcommand{\be}{\begin{equation}}
\newcommand{\ee}{\end{equation}}

\usepackage{lineno}

\begin{document}

%

 \singlespace
\title{On matrices whose exponential is a P-matrix
\thanks{The research of MM is
		supported in part 
		by a research grant 
		from  the Israel Science Foundation.  
}}

\author{Chengshuai Wu and Michael Margaliot\thanks{
		\IEEEcompsocthanksitem
	C. Wu is  with the School of Electrical  Engineering, Tel-Aviv University, Tel-Aviv~69978, Israel.
		\IEEEcompsocthanksitem
		M. Margaliot (Corresponding Author) is  with the School of Electrical  Engineering,
		and the Sagol School of Neuroscience, 
		Tel-Aviv University, Tel-Aviv~69978, Israel.
		E-mail: \texttt{michaelm@tauex.tau.ac.il}}}

\maketitle
 
\begin{center} 
								\today \\
\end{center}

\begin{abstract}
A matrix is called a P-matrix if all its principal minors are positive. P-matrices have found important applications in functional analysis, mathematical programming, and dynamical systems theory. We introduce  a new class of  real matrices   denoted~$\EP$. A   matrix is in~$\EP$ if and only if its matrix exponential is a P-matrix for all positive times. In other words, $A\in \EP$ if and only if  the transition matrix of the linear system~$\dot x=Ax$ is a P-matrix for any positive time~$t$. 
We analyze the properties of this new class of matrices  
and describe an application of our theoretical 
results to opinion dynamics. 
\end{abstract}



\section{Introduction}

A matrix~$A\in\R^{n\times n}$ is called a P-matrix if every principal minor of~$A$ is positive. In particular, the diagonal entries of~$A$ and the determinant of~$A$ are positive. 
The class of P-matrices, denoted~$\P$, includes 
  important matrix classes  such as  
positive-definite  matrices, non-singular 
M-matrices~\cite[Chapter~2]{topics_math_ana}, B-matrices~\cite{bmatrices},  totally positive matrices~\cite{pinkus,fallat2017total},
and diagonally dominant
matrices with positive diagonal entries. 

Fiedler and  Ptak~\cite{fiedler1962matrices} 
presented the first systematic study of  P-matrices. These matrices have found many applications in economics~\cite{nika_book}, dynamical systems~\cite{vol_comp},
and mathematical programming~\cite{murty_LCP}.
 For a   survey on P-matrices, see \cite[Ch. 4]{mat_posi_johnson}.

We briefly review some  of these applications.  
 P-matrices have been used to analyze  the injectivity of nonlinear mappings.
 Consider a $C^1$ mapping~$f: \Omega  \to \R^n$ with $\Omega \subseteq \R^n$. Let~$J(x):=\frac{\partial} {\partial x} f(x)$ denote the Jacobian of~$f$.
 It is natural to speculate that if~$\det(J(x))\not =0 $ for all~$x\in\R^n$ then~$f$ is injective. But this is not so in general.
  Gale and Nikaido~\cite{gale-nikaido1965}  proved that if~$J(x) \in \P$ for all~$x\in\Omega$ and~$\Omega$~is a rectangle then~$f$ is injective in~$\Omega$. 
 For generalizations of the Gale and Nikaido theorem, see e.g.,~\cite{gen_Gale_1,gen_Gale_2}.

Another important application of P-matrices is in the field of mathematical programming, in particular, the \emph{linear complementarity  problem} (LCP), which is a generalization of both linear programming and quadratic programming.  
For a  comprehensive treatment of the
LCP and, in particular, 
 its  numerous applications, see~\cite{LCP_cottle}. Some recent results on the relations between LCP and totally positive matrices are given in  \cite{choudhury2021}. 
Given~$M\in\R^{n\times n}$ and~$q\in\R^n$, the  $\text{LCP}(q,M)$ is: find (or conclude there is no)~$z \geq 0 $ such that
\begin{align}\label{eq:comp}
    & w:=M  z + q  \geq 0, \nonumber  \\
    & w^T z =0.
\end{align}
 Condition~\eqref{eq:comp}
 is called  the complementarity condition, since it implies that if~$z_i>0$ [$w_i>0$] for some index~$i$ 
 then~$w_i=0$ [$z_i=0$]. The LCP admits a unique solution for every~$q \in\R^n$ if and only if (iff)~$M \in\P$~\cite{murty_pmat}. 

P-matrices have also found applications in dynamical systems theory. We discuss their role in  a specific  model that
highlights an  intuitive interpretation
of~$P$-matrices. 
A fundamental  model in mathematical   ecology 
 is  the  
Lotka-Volterra  equations~\cite[Part~5]{sigmund_evolution_book}:
\be\label{eq:kol_sys}
\dot x_i = x_i(b_i + \sum_{j=1}^n a_{ij}x_j),\quad i=1,\dots,n. 
\ee
Here~$x_i(t)$ is the biomass of species~$i$ at time~$t$, $b_i $ is the growth rate of species~$i$, 
and~$a_{ij}   $ describes the interconnection between species~$j$ and species~$i$. It is clear that the non-negative orthant~$\R^n_+$ is an invariant set of the dynamics, and we assume throughout that~$x(0)\in\R^n_+$.  
Let~$A:=( a_{ij})_{i,j=1}^n$ and~$b:=\begin{bmatrix}
b_1&\dots& b_n
\end{bmatrix}^T$.
Any equilibrium point~$e$ of~\eqref{eq:kol_sys}  satisfies
\be\label{eq:eq_kol_sys}
 \diag(e_1,\dots,e_n) ( b+Ae)=0,
\ee
where~$\diag(c_1,\dots c_n) $ is the~$n\times n$ diagonal matrix with diagonal entries~$c_1,\dots,c_n$. 
The  Jacobian  of the vector field in~\eqref{eq:kol_sys}
is 
\begin{align}\label{eq:jlv}
J(x)&= 
\diag(g_1(x),\dots,g_n(x) ) +  \diag(x_1,\dots, x_n)A
 , 
\end{align}
where
\be\label{eq:defgi}
g_i(x):=b_i+\sum_{j=1}^n a_{ij}x_j.
\ee
If~$e$ is an equilibrium of~\eqref{eq:kol_sys}, and~$e_i=0$ for some index~$i$ then~\eqref{eq:eq_kol_sys}, \eqref{eq:jlv}, 
and~\eqref{eq:defgi}  imply 
that 
\[
(\zeta^i)^T J(e)   = g_i(e) (\zeta^i)^T,
\]
where~$\zeta^i$ is the $i$th canonical vector in~$\R^n$. Thus,  a necessary condition for the stability of~$e$ is that
\be\label{eq:lv_nec_stab}
g_i(e)=b_i+\sum_{j=1}^n a_{ij} e_j \leq 0, 
\ee
 for any $i$
  such that~$e_i=0$. Of course, \eqref{eq:lv_nec_stab} also holds (with an equality) for any~$i$ such that~$e_i>0$.  
Consider the 
$LCP (-b, -A)$, that is,  
find a vector~$y \geq 0$ 
such that~$ -A y  - b  \geq 0$ and~$ y^T (  -A y -b ) = 0$.  Then 
any stable 
equilibrium point~$e$ of~\eqref{eq:kol_sys} is a solution of this~$LCP$.  If~$(-A) \in \P$ then the LCP admits a unique solution for any~$b$, so~\eqref{eq:kol_sys}  has no more than a single stable equilibrium, for any~$b$.

Here, we 
introduce and analyze
a new  class of matrices: we  call a matrix~$A\in\R^{n\times n}$ an  \emph{exponential P-matrix} 
if~$\exp(At)$ is a~P-matrix for \emph{all}~$t \geq 0$. 
We denote this class of matrices by~$\EP$.
To the best of our knowledge, this type of matrices has not been studied  before. Our work is motivated in part 
by the work of 
Binyamin 
   Schwarz~\cite{schwarz1970} on totally positive differential  systems~(TPDSs). 
 Recall that a matrix  is called totally positive~(TP) if all its minors are positive~\cite{pinkus,total_book,gk_book}.  Schwarz considered the linear time-varying 
 system 
 \be \label{eq:linsys0}
 \dot x(t)=A(t)x(t), \quad x(t_0)=x_0,
 \ee
 where~$t\to A(t)$ is continuous.
 The solution of this system is~$x(t)=\Phi(t,t_0)x_0$
 where~$\Phi(t,t_0)$, the  transition matrix from~$t_0$ to~$t$,    is
 the solution at time~$t$ of the matrix differential equation
 \[
 \frac{d}{ds} \Phi(s,t_0)=A(s)\Phi(s,t_0),\quad \Phi(t_0,t_0)=I.
 \]
System~\eqref{eq:linsys0} is called a~TPDS if~$\Phi(t,t_0)$ is~TP for any pair~$t_0,t$ with~$t>t_0\geq 0$. 
Schwarz showed that if~\eqref{eq:linsys0} is a TPDS then  the variation diminishing property  of TP matrices~\cite{gk_book} can be used   to 
analyze the asymptotic behaviour of~$x(t)$. If~$A(t)\equiv A$ then~\eqref{eq:linsys0} is a TPDS iff~$A$ is Jacobi, that is,~$A$ is tri-diagonal with positive entries on the super- and sub-diagonals.
In other words,
\[
\exp(At) \text{ is TP for all } t>0 \text{ iff } A \text{ is Jacobi.}
\]
 It was recently shown~\cite{margaliot2019revisiting}  that TPDSs have important applications in the asymptotic analysis of the time-varying  nonlinear dynamical systems in the form~$\dot x=f(t,x)$ 
whose Jacobian~$J(t,x):=\frac{\partial}{\partial x} f(t,x)$ is a Jacobi matrix
 for all~$t,x$. See also~\cite{rami_osci} for the analysis of 
 discrete-time totally positive dynamical systems. The analysis of such systems builds on the use of the multiplicative- and additive-compounds of~$J$~\cite{comp_barshalom_omri}.

Note that if~$\Phi(t, t_0)$ is TP then in particular  it is a P-matrix. Thus, if~$A$ is Jacobi then~$A\in\EP$. 
However, the analysis  of matrices whose exponential is a P-matrix  
seems to be  more complicated than that of  matrices whose exponential is~TP due to the fact that~TP matrices are closed under multiplication, whereas~P-matrices are not. Another important difference   is that the class  of P-matrices is closed under matrix inversion, whereas the class  of TP matrices is not.

The contributions in this paper include the following. We introduce the  new class of matrices~$\EP$, provide various conditions  for a matrix to be in~$\EP$, describe  transformations that preserve~$\EP$,  describe an application to consensus systems, and analyze existing consensus algorithms in this new framework.

 The remainder of this paper is organized as follows. The next section reviews  some results  on P-matrices that  are used later on. 
 Section~\ref{sec:main} describes our main results. Section~\ref{sec:app} describes an application of our theoretical results  in the context of consensus algorithms~$\dot x=-Lx$, where~$L$ is a Laplacian matrix. 
 We  add a natural requirement on the dynamics of this system, called non-sign reversal, and show that it holds iff~$ L$ is an EP-matrix. 
The final section concludes and describes several directions for further research.

\emph{Notation}. Small [capital] letters denote column vectors [matrices]. For~$A,B\in \R^{n\times m}$, the
notation~$A\geq B$ [$A \gg B$] implies that~$a_{ij} \geq b_{ij}$ [$a_{ij} >  b_{ij}$] for all~$i,j$. We call~$A   \in \R^{n \times m}$ a non-negative [positive] matrix if $A \geq 0$ [$A \gg 0$]. The non-negative orthant in~$\R^n$ is~$\R^n_+:=\{x\in\R^n \st x\geq 0\}$.
 
 Consider a matrix $A \in \R^{n \times m}$ and fix an integer $k \in \{1, \dots, \min\{n,m\} \}$.
Let~$Q^{k,\ell}$ denote the set 
 of increasing sequences of~$k$ numbers from~$\{1, \dots, \ell \}$ ordered lexicographically. For example,
 \[
 Q^{2,3}=\{ \{1,2\}, \{1,3\}, \{2,3\} \}.
 \]
 A~$k$-minor of~$A$ is the determinant of some~$k \times k$
 submatrix of~$A$. 
 Each such submatrix  is defined by a set of row indices~$\alpha \in Q^{k,n}$ and column indices~$\beta \in Q^{k,m}$. This submatrix  
 is denoted by~$A[\alpha|\beta]$, and the corresponding minor is
 \[
 A(\alpha|\beta):=\det(A[\alpha|\beta]).
 \]
 In particular, $A(\alpha|\beta)$ is called a principle minor if $\alpha= \beta$, and a leading principle minor if $\alpha =\beta= \{1, 2,\dots, j \}$, for some~$ j \geq 1$. Similarly, $A[\alpha|\alpha]$ is called a principle submatrix of $A$. For $\alpha \in \Q^{k,n}$, let~$\bar \alpha:=\{1,\dots,n\}\setminus \alpha$ (we use set notation here, but we always  assume that  the elements in~$\bar \alpha$ are ordered in increasing order). 
 
The~$k$th \emph{multiplicative compound matrix} 
of~$A\in\R^{n\times m}$, denoted~$A^{(k)}$, is the~$\binom{n}{k}\times  \binom{m}{k}$ matrix that includes all  the   $k$-minors  ordered lexicographically. 
For example, if~$n=m=3$ and~$k=2$ then
\[
A^{(2)}= \begin{bmatrix}
A(\{12\}|\{12\}) & A(\{12\}|\{13\}) & A(\{12\}|\{23\}) \\ 
A(\{13\}|\{12\}) & A(\{13\}|\{13\}) & A(\{13\}|\{23\}) \\ 
A(\{23\}|\{12\}) & A(\{23\}|\{13\}) & A(\{23\}|\{23\}) 
\end{bmatrix}.
\]
Note that, 
by definition,~$A^{(1)}=A$ and if~$n=m$ then~$A^{(n)}=\det(A)$.
The $k$th \emph{additive compound matrix} of~$A \in \R^{n \times n}$
is    
\be \label{eq:poyrt}
	A^{[k]}:= \frac{d}{d \varepsilon}  (I+\varepsilon A)^{(k)} |_{\varepsilon=0}  .
\ee
Note that this implies that  $	A^{[k]}= \frac{d}{d \varepsilon} \left( \exp(A \varepsilon) \right)^{(k)}|_{\varepsilon=0}$.
There is an explicit formula for the entries of~$A^{[k]}$ in terms of the entries of~$A$ (see e.g. \cite{schwarz1970}). 
See~\cite{muldo1990,Hausdorff_contract} for more information  on  compound matrices and their applications to dynamical systems described by ordinary  differential equations. 
See also~\cite{comp_barshalom_omri}
for a recent tutorial on   applications of compound matrices  in systems  and control theory.  

A matrix $A \in \R^{n \times n}$ is called totally positive~(TP) [totally non-negative~(TN)] if~$A^{(k)} \gg 0 $  [$A^{(k)} \geq 0 $] for all~$k \in \{ 1, \dots, n\}$.  A matrix~$A\in \R^{n\times n}$ is called 
\emph{sign-symmetric} \cite{hersh2003} if for any~$k\in\{1,\dots,n\}$ and any~$\alpha,\beta \in Q^{k,n}$ we have
 $A(\alpha|\beta)A(\beta|\alpha)\geq 0$.  We say that~$A\in\R^{n\times n}$ is \emph{sign-pattern symmetric} if $a_{ij}a_{ji} \geq 0$ for all $i,j \in \{1, \dots, n\}$. Then~$A$ is sign-symmetric iff~$A^{(k)}$ is sign-pattern symmetric for all~$k \in \{1, \dots, n-1\}.$
  
 \section{Preliminaries: P-matrices}
A matrix~$A\in\R^{n\times n}$ is called a P-matrix if all its principal minors are positive, that is,   $A^{(k)}$ has positive diagonal entries for all $k \in \{1, \dots, n\}$. 
We list some transformations that preserve the P-matrix property. Recall that~$S\in\R^{n\times n}$ is called a signature matrix if it is a diagonal matrix, and  every diagonal entry is either one or minus one. Then~$S^{-1}=S$.
A matrix~$V\in \R^{n \times n}$ is called a permutation matrix if it has exactly one entry of 1 in each row and each column, and zeros elsewhere. Then~$V^{-1} = V^{T}$.
\begin{Theorem}\cite{gen_p_matrices} \label{thm:sign}
Suppose that~$A\in  \R^{n\times n}$ is a P-matrix. Then
\begin{enumerate}
    \item $A^T\in \P$. \label{item:at}
    \item  \label{item:qpq} If~$V$ is a permutation matrix then~$V A V^T \in \P$.
\item \label{item:d1d2}
If~$D_1,D_2$ are diagonal matrices and~$D_1D_2$ has positive diagonal entries   then~$D_1 A D_2 \in \P$.
 \item \label{item:sas} 
    If~$S$ is a signature  matrix then~$S   A S \in \P$.
\item \label{item:aplusd} If~$D $ is a non-negative diagonal matrix    then~$A+D \in \P$.
\item \label{item:homotopy}
If~$D$ is a 
diagonal matrix with~$0\leq D \leq I$ then~$D  +(I-D)  A \in \P$.
\item \label{item:inverse}
$A^{-1} \in \P.$
\end{enumerate}
\end{Theorem}

Since~$\exp(At)=I+\sum_{k=1}^\infty \frac{A^k t^k}{k!}$,  $\exp(A t )\in\P$ for all~$t>0$ sufficiently small. The next result  provides more information on the value~$t$ for which~$\exp(At) \in \P$ when~$A$ has non-negative entries. 

\begin{Theorem}\cite{gen_p_matrices}
If~$A\geq 0$  then~$\exp(At)\in\P$ 
for all~$t\in[0,\frac{1}{\rho(A)})$, where $\rho(A)$ is the spectral radius of~$A$. 
\end{Theorem}

The next result summarizes several  necessary and sufficient conditions for a matrix to be a $P$-matrix. 

 \begin{Theorem} \cite{fiedler1962matrices,MOYLAN197753} \label{thm:pmat} 
 Let~$A \in\R^{n\times n}$. The following six  conditions are equivalent. 
 \begin{enumerate}
     \item \label{cond:allmposi}
     Any   principal minor  of~$A$ is positive, that is,~$A\in\P$;
     \item  \label{cond:signrev}
     For any~$x\in \R^n\setminus\{0\}$ there exists an index~$k\in\{1,\dots,n\}$ such that~$x_k(Ax)_k>0$;
     \item \label{cond:diagsrev} For any~$x\in \R^n\setminus\{0\}$ there exists a positive diagonal matrix~$D(x)$ such   such that~$x^T  D(x) A x>0$;
     \item \label{cond:diagsrevnn}
     For any~$x\in \R^n\setminus\{0\}$ there exists a non-negative diagonal matrix~$H(x)$ such   such that~$x^T  H(x) A x>0$;
 \item \label{cond:allp} Every real eigenvalue of~$A$ and of any principal submatrix
 of~$A$
 is positive;
  \item \label{cond:aust}  For any signature matrix~$S$ there exists a vector~$x\gg 0$ such that~$SAS x \gg 0$.
 \end{enumerate}
 \end{Theorem}

 \begin{Remark} \label{re:signrev}
 For~$A\in\R^{n\times n}$, let
 \[
 \text{rev}(A):= \left \{x\in\R^n: x_i(Ax)_i \leq 0 \text{ for all } i\in\{1,\dots,n\}\right  \}.
 \]
 Roughly speaking, this is the set of vectors that~$A$ ``sign reverses''.
Condition~\ref{cond:signrev}) in Theorem~\ref{thm:pmat} implies that~$A$ is a P-matrix iff~$\text{rev}(A)=\{0\}$. 
This is known as the  sign non-reversal property (see, e.g.,~\cite{Apoorva2021}).
\end{Remark}

\begin{Remark}
Recall  that~$A\in\R^{n\times n}$ is called stable (or Hurwitz) if~$\real(\lambda)<0$ for any eigenvalue~$\lambda$ of~$A$. It is called diagonally stable if there exists a positive diagonal matrix~$D$ such that~$DA+A^TD  $ is negative-definite~\cite{diag_stab_book}. 
Condition~\ref{cond:diagsrev}) in Theorem~\ref{thm:pmat} implies that if $A  $ is diagonally stable   then $(-A) \in \P$. 
Recall~\cite{posi-tutorial}  that if~$A$ is Hurwitz and Metzler   then it is diagonally stable, so~$(-A)\in\P$
(see~\cite{box_inv2009} for an application of this property to determining box invariance of dynamical systems). The class of matrices
\[
\{A\in\R^{n\times n}:-A\in\P\}
\]
is sometimes denoted by~$\P^{(-)}$ (see~\cite{banaji2007}).
\end{Remark}

Recall that~$A\in\R^{n\times n}$ is called a Q matrix if for any~$k\in\{1,\dots,n\}$ 
the sum of all the~$k$ principal minors of~$A$ is positive. In particular, any P matrix is a Q matrix. 
The next result provides
a  necessary spectral condition for~$A \in \R^{n \times n}$ to be a Q-matrix. For~$z\in \C$, let 
$\arg(z) \in (-\pi,\pi]$ denote the argument of~$z$. 
\begin{Theorem}\cite{kellogg1972complex} \label{thm:arg}
Let $A \in \R^{n \times n}$, with~$n\geq 2$,
be a Q matrix (and thus a P matrix). Then any eigenvalue~$\lambda$ of~$A$  satisfies
\be \label{eq:nec_pmat}
| \arg(\lambda ) | < \pi - \frac{\pi}{n} .
\ee
\end{Theorem}
In other words, there is a  wedge
around the negative $x$-axis which is free from eigenvalues of~$A$.
In particular, a P-matrix cannot have 
a real and negative  eigenvalue (the latter is   also implied by Condition~\ref{cond:signrev}) in Theorem~\ref{thm:pmat}).
 
  It is well-known that the  problem of testing whether a given~$n\times n$ matrix is   a
P-matrix, formulated as a decision problem, is co-NP-complete in~$n$~\cite{coxson}.

The next section describes our main results. 
 \section{Main Results}\label{sec:main}
We begin by introducing a new 
class of matrices that we call
\emph{exponential P-matrices}. 
 \begin{Definition}\label{def:epmat}
 A matrix~$A\in\R^{n\times n}$ is called an EP-matrix if
\[
 \exp(At) \in \P \text{ for all } t\geq 0.
\]
 \end{Definition}
 Let~$ \EP$ denote the class of exponential~$P$ matrices.

One motivation for the notion of EP-matrices is a dynamical systems generalization  of the  sign-reversal property described in 
Remark~\ref{re:signrev}.
Indeed, if~$\dot x= Ax $ then the following two properties are equivalent:
\begin{enumerate}
    \item $A$ is an EP-matrix;
    \item for any~$x(0)\in\R^n\setminus\{0\}$ and any time~$t \geq 0 $ there exists at 
    least one index~$i=i(t,x(0))$ such that
    \[
    x_i(0)x_i(t)>0.
    \]
    
\end{enumerate}
In other words, at each time~$t$ there is at least one state-variable 
that has the same sign at time~$t $ as in time~$0$. 
We will describe an application of this property  to   opinion dynamics  in Section~\ref{sec:app}.

  Definition~\ref{def:epmat} only considers  non-negative times. However, using Condition~\ref{item:inverse}) in Theorem~\ref{thm:sign} yields the following result.
 \begin{Corollary}
 A matrix~$A\in\R^{n\times n}$ is~an EP-matrix iff
 \be\label{eq:all_time}
 \exp(At) \in \P \text{ for all } t\in\R.
 \ee
 \end{Corollary}
 \begin{IEEEproof}
 If~\eqref{eq:all_time} holds then clearly~$A \in \EP$. Now suppose that~$A \in \EP$, that is,~$\exp(At)\in\P$ for any~$t\geq 0$. 
 Fix~$s<0$. Then
$
 \exp(sA)= (\exp(-sA))^{-1}, 
 $ and
 since~$(-s)>0$ and the inverse of a P-matrix is a P-matrix, we conclude  that~$\exp(sA)\in \P$.
 \end{IEEEproof}
 
 \begin{Example}\label{exa:minus_jacobi}
 Consider $A=\begin{bmatrix}
 0& -1 &0 \\
-2& 0 &-2\\
 0&-1&0
 \end{bmatrix} $. Note that~$(-A)$ is Jacobi. A calculation gives
 \[
 \exp(At)= \begin{bmatrix}
 \cosh^2(t) & -\sinh(t)\cosh(t)  &\sinh^2(t) \\
-\sinh(2t)& \cosh(2t) &-\sinh(2t)\\
 \sinh^2(t)&-\sinh(t)\cosh(t)&\cosh^2(t)
 \end{bmatrix}, 
 \]
and 
 \[ 
 ( \exp(At))^{(2)} = 
 \begin{bmatrix} 
 \cosh^2(t) & -\sinh(2t)  &\sinh^2(t) \\
-\sinh(t)\cosh(t) & \cosh(2t) &-\sinh(t)\cosh(t)\\
 \sinh^2(t)&-\sinh(2t) &\cosh^2(t)  
 \end{bmatrix}.
 \]
 The diagonal entries of~$(\exp(At))^{(k)}$ are the principal minors of order~$k$ of~$\exp(At)$. 
 We see that the  principal minors  of order one and two of~$\exp(At)$ are positive for all~$t\geq 0$,
 so~$A\in\EP$. Note that the principal minors are   positive for any~$t\leq 0$ as well.
However,    some non-principal  minors of~$\exp(At)$ may take  negative values, as~$A$ is not Jacobi.
 \end{Example}
 
 
 More generally, since~$\exp(-At)$ is the inverse of~$\exp(At)$,
 it follows from   Jacobi's identity~\cite[Chapter~0]{matrx_ana} that for any~$k\in\{1,\dots,n-1\}$~and any~$\alpha\in Q^{k,n}$, we have
 \[
 (\exp(-At))(\alpha|\alpha)=
 \frac{ (\exp(At)) (\bar \alpha|\bar \alpha) }{\det(\exp(At))},
 \]
where~$\bar\alpha:=\{1,\dots,n\}\setminus\alpha$.
This shows, in particular, that~$A\in \EP$ iff~$-A \in \EP$.

The next result lists  
  some transformations that preserve the~EP-property.  
\begin{Theorem} \label{thm:epsign}
Suppose that~$A\in \R^{n\times n}$ is~an EP-matrix. Then
\begin{enumerate}
\item $c A \in \EP$ for any~$c\in\R$. \label{item:scaleda}
\item $A^T\in \EP$. \label{item:at4}
\item  \label{item:qpq4} If~$Q$ is a permutation matrix then~$QAQ^T \in \EP$.
\item \label{item:dadm1} 
    If~$D$ is a positive diagonal matrix    then~$D   A D^{-1}\in \EP$.
\item \label{item:sas4} 
    If~$S$ is a signature  matrix then~$S   A S \in \EP$.
\item \label{item:aplusd4}
If~$D $ is a  diagonal matrix  and~$DA=AD$  then~$A+D \in \EP$.
\end{enumerate}
\end{Theorem}

 \begin{IEEEproof}
To prove \ref{item:scaleda}), fix~$c\in \R$. 
Then~$\exp ((c A) t) =\exp(A (ct))$, and since~$A\in \EP$, $ \exp(A (ct))\in \P$.
The proof of~\ref{item:at4}) follows from the fact that~$ \exp(A^T t)=(\exp(At))^T $.
 To prove    \ref{item:qpq4}), note that if~$Q$ is  a permutation matrix then
 \[
    \exp(QAQ^T t)=Q \exp(At) Q^T,
 \]
 and since~$\exp(At)\in\P$, this implies that~$QAQ^T \in \EP$. The proofs of~\ref{item:dadm1}) and~\ref{item:sas4}) are similar to the proof of~\ref{item:qpq4}).
 To prove~\ref{item:aplusd4}),  let~$B:=A+D$. Then~$\exp(Bt)=\exp(Dt)\exp(At)$. Since~$\exp(Dt)$ is a positive diagonal matrix, and~$\exp(At)\in \P$, we conclude from  
 Condition~\ref{item:d1d2}) in Theorem~\ref{thm:sign}
 that~$\exp(Bt) \in \P$.
  \end{IEEEproof}
 
The following example shows that 
 the EP property is in general not invariant under   similarity transformations. 
\begin{Example}\label{exa:ninv}
The matrix~$
 A=\begin{bmatrix}
 0&1\\0&0
 \end{bmatrix}
 $
 is EP, as~$\exp(At)=\begin{bmatrix}
 1&t\\0&1
 \end{bmatrix}$,
 but for~$T=\begin{bmatrix}
 -1&-1\\1&0
 \end{bmatrix}$
 we have
 \[
 B:=TAT^{-1} = \begin{bmatrix}
 1&1\\-1&-1
 \end{bmatrix}
\]
and $B$ is not EP, as
\[
\exp(Bt)=T\exp(At) T^{-1} =\begin{bmatrix}
1 + t& t  \\  -t& 1 - t
\end{bmatrix}
\]
and for any~$t\geq 1$ this matrix has a non-positive  diagonal entry, so it is not a P~matrix.
\end{Example}

 \subsection{Some classes of EP-matrices}
 In general, it seems difficult to  determine if a given matrix is~$\EP$. However, we show that some specific  classes of matrices   are $\EP$. 
 We first study the case of  $2 \times 2$ matrices  where it is possible to give a necessary and sufficient  sign-pattern condition for~$\EP$. 

\begin{Proposition}\label{prop:case2by2}
  Suppose that~$A\in\R^{2\times 2}$.
  Then~$A\in\EP$   iff~$a_{12}a_{21}\geq 0$.
\end{Proposition}
In other words,~$A\in\EP$ iff~$A$ has the sign pattern
\be\label{eq:sp}
\begin{bmatrix}
* &\geq 0\\
\geq 0 &* 
\end{bmatrix}
\text{ or } 
\begin{bmatrix}
* &\leq 0\\
\leq 0 &* 
\end{bmatrix},
\ee
where~$*$ denotes ``don't care'', that is, iff~$A\in\R^{2\times 2}$ sign-pattern symmetric. Note that the first [second]
sign pattern here corresponds to~the linear dynamical system $\dot x=Ax $ being a cooperative [competitive]
system~\cite{hlsmith}.

\begin{IEEEproof}
If~$a_{12} a_{21}=0$ then~$A$ is triangular and so is $\exp(At)$. In this case, the diagonal entries of $\exp(At)$ are $\exp(a_{ii}t)$, $i =1, 2$. Hence,~$A\in \EP$. 

If~$a_{12}>0$ and~$a_{21}>0$
then~$A$ is Jacobi, so~$A\in\EP$. If~$a_{12}<0$ and~$a_{21}<0$, then $-A \in \EP$ and thus~$A \in \EP$. 

To complete the proof,  we need to show that if
\be\label{eq:condnreg}
a_{12}a_{21}<0
\ee
then~$A\not \in\EP$. 

Let~$B(t):=\exp(At)$.
Since~$A\in\R^{2\times 2}$, $B(t)$ is a P-matrix iff its diagonal  entries~$b_{11}(t)$
and~$b_{22}(t)$   are positive. 
Let~$s:=\trace(A)/2$, and write
$
A=s I +(A-s I) .
$
Note that $\trace(A-sI) = 0$.
Then
\[
B(t)= \exp(s t) \exp((A-s I)t).
\]
The term~$\exp(st)$ has no effect on the signs of~$b_{11}(t)$, $b_{22}(t)$, so we may assume that
\[
A=\begin{bmatrix} a_{11} &a_{12}\\
a_{21} &- a_{11}
\end{bmatrix}.
\]
The eigenvalues of~$A$ are
$
\pm \sqrt{a_{11}^2+a_{12}a_{21}}.
$
Assume that~\eqref{eq:condnreg} holds. 
We consider three cases. 

{\sl Case 1.} Suppose that~$ a_{12}a_{21}<-a_{11}^2$. In this case, the eigenvalues
  of~$A$ are purely imaginary and Theorem~\ref{thm:complex_eig} implies that~$A\not \in \EP$. 

 {\sl Case 2.} Suppose that~$a_{12}a_{21}=-a_{11}^2$. Then~$\det(A)=\trace(A)=0$, so~$A^2=0$ and thus
\[
B(t)=I+At=\begin{bmatrix}
1+a_{ 11} t& a_{12} t\\
a_{21} t& 1-a_{11} t
\end{bmatrix} . 
\]
Since~$a_{12}a_{21}<0$,~$a_{11}\not =0$, and   we conclude that~$B(t)\not\in \P$   for any~$t\geq 0$ sufficiently large, so~$A\not\in \EP$.

{\sl Case 3.} Suppose that~$ a_{12}a_{21}>-a_{11}^2$. Then~$\det(A)=-a_{11}^2-a_{12}a_{21} <0$, and the diagonal entries of~$B(t)$ are
\begin{align*}
  b_{11}(t)  &=\cosh(s t  ) \left(1+  \frac{ a_{11 }}{ s}  \tanh (s t)\right ),\\ 
b_{22}(t) &= \cosh(s t  )\left(1 -  \frac{ a_{11 }}{ s}  \tanh (s t)\right ) ,
\end{align*}
where~$s:=\sqrt{-\det(A ) }$.
Since~$a_{12}a_{21}<0$,  $\frac{ |a_{11 }|}{ s}>1 $, so~$B(t) \not \in \P $ for all~$t> 0$ sufficiently large, so~$A\not \in \EP$. 
  \end{IEEEproof}
 
 \begin{Example}
 Consider the matrix~$A(  w):=\begin{bmatrix}
 -1&w\\-w&-1
 \end{bmatrix}$, with~$ w>0$. 
 Then
 \[
 \exp(At)=\exp(- t) \begin{bmatrix}
 \cos(wt) &\sin(wt)\\
 -\sin(wt)& \cos(wt)
 \end{bmatrix},
 \]
 so~$\exp(At)\in \P$ for~$t\in[0,\frac{\pi}{2 w})$, but
 \begin{align*}
     \exp(A  \frac{\pi}{2 w})= \exp(- \frac{\pi}{2w}) \begin{bmatrix}
0&1\\
 -1& 0
 \end{bmatrix} \not \in \P,
     \end{align*}
    so~$A\not\in \EP$. Note that when~$w \to 0$,  the $L_2$ norm of~$A$ converges to one, yet the first  time where we ``loose'' the~$\EP$ property, namely,~$\frac{\pi}{2 w}$ converges to infinity. 
 \end{Example}

The next result    describes several classes of matrices that are~$\EP$.
\begin{Theorem} \label{thm:epmc}
Let $A \in \R^{n \times n}$. Any one of the following three conditions implies that~$A\in\EP$. 
\begin{enumerate}[1)]
    \item  \label{item:tria}  
 $A$ is    lower  triangular or upper triangular.
    \item \label{item:sym}  $A$ is symmetric.
    \item \label{item:jaco}   $A$ is Jacobi.
\end{enumerate}
\end{Theorem}

\begin{IEEEproof}
Pick an arbitrary~$t\geq 0$ and let~$B:=\exp(At)$. We will show that if any one of the conditions in the  theorem holds then~$B\in \P$.

To prove~\ref{item:tria}), note that the (lower or upper) triangular structure is preserved under summation and multiplication. Hence, $B = \sum_{\ell=0}^\infty \frac{A^\ell t^\ell}{\ell !}$ is also a triangular matrix. Since the eigenvalues of a triangular matrix are its diagonal entries, the diagonal entries of $B$ are $\exp(a_{ii}t)$,  $i = 1, \dots, n$. Pick~$k \in  \{1, \dots, n\}$ and $\alpha \in Q^{k,n}$. Then 
$
B(\alpha | \alpha) = \prod_{i \in \alpha} \exp(a_{ii}t)>0,
$
so~$B\in \P$.
  
To prove~\ref{item:sym}),
 suppose that $A$ is symmetric.  Then all the eigenvalues~$\lambda_i$ of~$A$ are real, $B$ is symmetric and  the eigenvalues of~$B$ are $\exp(\lambda_i t)>0$, $i=1,\dots,n$. Thus,~$B$ is positive-definite, and by Theorem~\ref{thm:pmat},  $B \in \P$.

To prove~\ref{item:jaco}), suppose that~$A$ is Jacobi. Then~$\exp(At)$ is~TP for all~$t>0$, and in particular~$A$ is~$\EP$. Here, we also give a more direct  proof. Since~$A$ is Jacobi, it has the form 
\[
A  =
\begin{bmatrix}
a_1 & b_1 & 0 & 0 & \cdots & 0 & 0 & 0 \\
c_1 & a_2 & b_2 & 0 & \cdots & 0 & 0 & 0 \\
0 & c_2 & a_3 & b_3 & \cdots & 0 & 0 & 0 \\
\vdots & &  &  & \vdots  &  & & \vdots \\
0 & 0 & 0 & 0 & \cdots & c_{n-2} & a_{n-1} & b_{n-1} \\
0 & 0 & 0 & 0 & \cdots & 0 & c_{n-1} & a_n 
\end{bmatrix},
\]
with~$b_i, c_i > 0$, $i = 1, \dots, n-1$. Consider  the positive diagonal matrix~$D:= \diag(1, \sqrt{\frac{b_1}{c_1}}, \dots, \sqrt{\frac{b_1 b_2 \cdots b_{n-1}}{c_1 c_2 \cdots c_{n-1}}})$. Then~$D A D^{-1}$ is symmetric. Hence, $\exp(DA D^{-1}t) = D \exp(At) D^{-1} \in \P.$  Theorem~\ref{thm:sign} implies that $\exp(At) \in \P$. 
\end{IEEEproof}
The three classes of matrices in Theorem~\ref{thm:epmc}   all have real eigenvalues. One may conjecture that if all the eigenvalues of~$A \in \R^{n \times n}$ 
are real then $A \in \EP$.  However, the matrix~$B$ in Example~\ref{exa:ninv} above
shows that this conjecture is false.  
 
The next result shows  that having all real eigenvalues is   a necessary condition for EP-matrices.
\begin{Theorem}\label{thm:complex_eig}
 If $A \in \R^{n \times n}$ has  a complex (non-real) eigenvalue  then $A \notin \EP$. 
\end{Theorem}

\begin{IEEEproof}
Let $\lambda  =  a + j b$, where $b \neq 0$ and $j :=  \sqrt{-1}$, be a complex  eigenvalue of $A$. Since the complex eigenvalues of a real matrix always occur in conjugate pairs, we can assume that~$b>0$. Now, $\exp(At)$ has the eigenvalue 
\[
\exp(\lambda t) = \exp(at) \exp(jbt) .
\]
For $ t = \frac{(n-1)\pi}{n b } > 0$, we have $\arg(\exp(\lambda t)) = \pi - \frac{\pi}{n}$, and   Theorem~\ref{thm:arg} implies that~$\exp(A t )\not\in \P $, so~$A \notin \EP$.
\end{IEEEproof}

\begin{Remark}
An immediate implication of Theorem~\ref{thm:complex_eig} is that if~$A\in\R^{n\times n}$ is skew-symmetric then it is not an EP-matrix. If~$A \in \EP$, then Theorem~\ref{thm:complex_eig} also implies that the linear dynamical system $\dot x = Ax$ does not have limit cycles. 
\end{Remark}

The following result  
  provides a necessary and sufficient condition for~EP in terms of the diagonal entries of  additive  compound matrices. Recall that if~$A\in\R^{n\times n}$ then~$A^{[k]} \in \R^{r\times r}$, with~$r:=\binom{n}{k}$.
\begin{Proposition}\label{prop:necsuffk}
Let~$A\in\R^{n\times n}$. Then~$A\in \EP$
iff for any~$t\geq 0$ and any~$k\in\{1,\dots,n-1\}$, we have
\[
(\exp (A^{[k]} t)  )_{ii} >0, 
\text { for all } i\in\{1,\dots,\binom{n}{k}\}.
\]
\end{Proposition}
\begin{IEEEproof}
Recall (see e.g.~\cite{muldo1990}) that
\[
(\exp( A  ))^{(k)} = \exp( A^{[k]} )  .
\]
 Thus, if~$\alpha$ is the~$i$th element 
in~$ Q^{k,n}$ then the corresponding principal minor of~$\exp(At)$ satisfies 
\[
(\exp( A  ))(\alpha|\alpha)= (\exp( A^{[k]} )  )_{ii} ,
\]
and this completes the proof.
\end{IEEEproof}

The next result provides a sufficient  condition for~EP
based on sign-symmetry of~$\exp(At)$.
\begin{Proposition} \label{prop:signsym}
Suppose that~$A\in\R^{n\times n}$ has only real eigenvalues and that
\be\label{eq:signsymmc}
\exp(At) \text{ is sign-symmetric for all } t\geq 0.
\ee
Then~$A\in\EP$.
\end{Proposition}
Note that by
Theorem~\ref{thm:complex_eig},
the assumption that~$A$ has only real eigenvalues is a necessary condition for~$\EP$.
\begin{IEEEproof}
Fix~$t\geq0$. 
Since~$A$ has real eigenvalues, $\exp(At)$
has real and positive eigenvalues.
In particular~$\exp(At)$ is positively stable (i.e.,~$(-\exp(At))$ is a Hurwitz matrix). Recall that  for the class of
sign-symmetric matrices, positivity of principal minors and positive stability are equivalent (see, e.g.,~\cite{hersh2003}), so we conclude that~$\exp(At)\in\P$.
\end{IEEEproof}

\begin{Remark}
Fix~$k\in\{1,\dots,n\}$. It follows from~\eqref{eq:poyrt} that 
\[
( \exp(A t) )^{(k)} = I_r + t A^{[k]} +o(t), 
\]
with~$r:=\binom{n}{k}$,
so if condition~\eqref{eq:signsymmc} holds
then~$I_r+t A^{[k]}$ is sign-pattern symmetric for any $ t>0$ sufficiently small, that is, 
$  A^{[k]}$ is sign-pattern symmetric.
\end{Remark}

Prop.~\ref{prop:signsym} suggests the following question: suppose that~$A\in\R^{n\times n}$ has only real eigenvalues and is sign-symmetric. Is~$\exp(At)$ sign-symmetric for all~$t\geq 0$? The following example shows that in general the answer is no.

\begin{Example}
Consider the matrix~$A=\begin{bmatrix}
-1/2& 1&1 \\
1&10&1\\
1&0&20
\end{bmatrix}$. It is straightforward to verify that: $A$ has three real eigenvalues, $A$ is sign-symmetric, but~$
\exp(0.1 A) 
$ is not 
sign-symmetric.
\end{Example}

There are however special cases where sign-symmetry of~$A$ does allow us to apply Prop.~\ref{prop:signsym}.   
Recall that a matrix is called~P$_0$ if all its principal minors are non-negative (see, e.g.,~\cite{hersh2003}). 

\begin{Corollary}
Suppose that~$A\in\R^{n\times n}$ is sign-symmetric, and is a projection matrix (i.e.,~$A^2=A$). Then~$A\in\EP$. 
\end{Corollary}

Note that we do not require here explicitly that~$A$ has only real eigenvalues. Indeed, 
since~$A$ is a projection matrix, all its  eigenvalue     either one or zero.
\begin{IEEEproof}
Fix~$t\geq0 $, and let~$B:=\exp(At) - I$. Since~$A^2=A$,
\begin{align*}
   B&= \sum_{k=1}^\infty \frac{A t^k}{k!}\\
                &= ( \exp(t) - 1 )A.
\end{align*}
Since $A$ is sign-symmetric, $A^2$ is a P$_0$-matrix \cite[Lemma~2.1]{hersh2003}, and since~$A^2=A$, $A$ is a P$_0$-matrix.
Hence, $B$ is also a P$_0$-matrix.
Now~\cite[Theorem~1.5]{FIEDLER1966} implies that~$\exp(At) = I + B$ is a P-matrix.
\end{IEEEproof}

\begin{Example}
Consider the matrix $A =(1/4) \begin{bmatrix}   4 & -3& -3 \\ 0 & 1 & -3 \\ 0 & -1 & 3  \end{bmatrix}.$ This is a sign-symmetric projection matrix.
A calculation gives
\[
\exp(At)=(1/4)\begin{bmatrix}
4\exp(t)& 3 - 3\exp(t)&  3 - 3\exp(t)\\
0&   3+\exp(t) & 3 - 3\exp(t)\\
0&   1 - \exp(t)& 1+3\exp(t)
\end{bmatrix},
\]
and
\[
(\exp(At))^{(2)}=(\exp(t)/4)\begin{bmatrix}
3+\exp(t)&*&*\\
*& 1+3\exp(t) &*\\
*&*& 4
\end{bmatrix},
\]
where~$'*'$ denotes values that are not relevant. 
Since all the diagonal entries of~$\exp(At)$ and~$(\exp(At))^{(2)}$
are positive for all~$t$, 
 $A \in \EP$.
\end{Example}

 \subsection{Generalizations of Theorem~\ref{thm:epmc}}

This section provides some generalizations of Theorem~\ref{thm:epmc}. Specifically, it is straightforward to verify that the three classes of matrices in Theorem~\ref{thm:epmc} are   special cases of Theorem~\ref{prop:reducible}, Theorem~\ref{prop:diag}, and Theorem~\ref{thm:genjacobi}, respectively. Furthermore, all the matrix classes in Theorem~\ref{thm:epmc} are included in Theorem~\ref{thm:signsym}. 

\begin{Theorem}
\label{prop:reducible}
Suppose that
\be\label{eq:reducform}
A=\begin{bmatrix}
B &C \\ 0& D 
\end{bmatrix},
\ee
where~$B\in\R^{n\times n}$,
$C \in \R^{n\times m}$,
$D \in \R^{m\times m}$,
and~$0$ denotes an~$m\times n$ matrix of zeros. Then~$A\in \EP$ iff~$B,D\in\EP$.
\end{Theorem}
\begin{IEEEproof}
First, we show~$A \in \EP$ implies~$B,D\in\EP$. Note that 
\[
\exp(At) = 
\begin{bmatrix}
\exp(Bt) & * \\
0 & \exp(Dt) 
\end{bmatrix},
\]
where $*$ denotes ``don't care''. Since $\exp(At)$ is a P-matrix, every   principal submatrix
of~$\exp(At)$ is also a P-matrix. Hence, $B,D\in\EP$.

To prove the converse implication, assume that~$B,D \in \EP$. 
Fix~$t\geq 0$ and~$v\in \R^{n+m }\setminus\{0\}$.
Let~$w:=\exp(At)v$, i.e. the solution at time~$t$ of~$\dot x =Ax$, $x(0)=v$. 
By Theorem~\ref{thm:pmat},
it is enough to show that there exists a~$k\in\{1,\dots,n+m\}$ such that~$w_k v_k>0$.
We consider two cases. 

\noindent \emph{Case 1.}
Suppose that at least one  of~$v_{n+1},\dots,v_{n+m}$ is non zero.
Then 
\be \label{eq:case1}
\begin{bmatrix}
w_{n+1}\\ \vdots\\ w_{n+m}
\end{bmatrix}=\exp(Dt)\begin{bmatrix} v_{n+1}\\ \vdots\\ v_{n+m}
\end{bmatrix},
\ee
and the fact that~$D\in\EP$ implies that there exists~$i\in\{n+1,\dots,n+m\}$ such that $w_i v_i >0$.

\noindent \emph{Case 2.}
Suppose that~$v_{n+1}=\dots=v_{n+m}=0$. Then
\be \label{eq:case2}
\begin{bmatrix}
w_{1}\\ \vdots\\ w_{n}
\end{bmatrix}=\exp( B t)\begin{bmatrix} v_{1}\\ \vdots\\ v_{n}
\end{bmatrix},
\ee
and since~$v\not =0$ and~$B\in\EP$, there exists~$i\in\{1,\dots,n\}$ such that $w_i v_i >0$.
Hence,~$A \in \EP$.
\end{IEEEproof}

\begin{Remark}
Any reducible matrix is similar to~\eqref{eq:reducform} up to a  permutation similarity. Since the EP-property is preserved under permutation similarity (see Theorem~\ref{thm:epsign}), Proposition~\ref{prop:reducible} shows that when analyzing EP-matrices it is enough to study irreducible matrices.
\end{Remark}

\begin{Example}
Let
\be\label{eq:adoesnt}
A = 
\begin{bmatrix}
-4 & 1 & -6 \\
2 & -5 & 6 \\
0 & 0 & -9
\end{bmatrix}.
\ee
Note that $A$ does not belong to any of the  three matrix
classes   in Theorem~\ref{thm:epmc}, but it  has the form~\eqref{eq:reducform} with~$B=\begin{bmatrix}
-4 & 1  \\
2 & -5 
\end{bmatrix}$ and~$D=-9$. Since~$B,D \in \EP$, so is~$A$. Indeed, a calculation gives that 
\begin{align*}
3\exp(At) = & \left [
  \begin{smallmatrix}
 2 \exp(-3t) + \exp(-6t)  & * & * \\
*&  \exp(-3t) +  2 \exp(-6t)   & * \\
* & * & 3\exp(-9t)
\end{smallmatrix} \right ],
\end{align*} 
and that~$3(\exp(At))^{(2)}$ is equal to
\begin{align*}
\left [ \begin{smallmatrix}
3\exp(-9t) & * & * \\
*&  2 \exp(-12t) + \exp(-15 t)  & * \\
* & * &  \exp(-12t) + 2 \exp(-15 t) 
\end{smallmatrix}\right ] ,
\end{align*}
so all the principal minors of~$\exp(At)$
are positive for all~$t$.
\end{Example}

The next result provides a sufficient condition for~$A=TDT^{-1}$, with~$D$ a diagonal matrix, to be~$\EP$. This condition is a kind of symmetry condition for the minors of~$T$.  For~$\alpha ,\beta  \in Q^{k,n}$, let
\[
s(\alpha,\beta):=(-1)^{\sum_{j\in\alpha}  j+\sum_{i\in\beta} i }  ,
\]
that is, the signature of~$\alpha+\beta$.

\begin{Theorem}\label{prop:diag}
Suppose that~$A\in\R^{n\times n}$  has real eigenvalues~$\lambda_1,\dots,\lambda_n$ and is diagonalizable, that is, there exists a non-singular matrix $T\in\R^{n\times n}$ such that $A=T D T^{-1}$ with $D :=\diag(\lambda_1,\dots,\lambda_n)$.
If for any~$k\in\{1,\dots, \lceil \frac{n-1}{2} \rceil \}$ and any~$\alpha  \in Q^{k,n}$ we have
\be\label{eq:tadjugate}
\det(T)
s(\alpha,\beta)    T(\alpha|\beta)   
 T(\bar\alpha|\bar\beta)\geq 0 \text{ for all } \beta\in Q^{k,n} ,
\ee
and this holds with an inequality  for at least one~$\beta\in Q^{k,n}$, 
then~$A\in\EP$.
\end{Theorem}

\begin{Example} \label{exa:Tdiag}
Let~$T=\begin{bmatrix}
3 & 3 & -3\\
3  &-1&-1\\
1 & 1 &1
\end{bmatrix}$.
Then~$\det(T)= -24 <0$.
Calculating
\[
r(\alpha,\beta)
:=s(\alpha,\beta)  T(\alpha|\beta)   
 T(\bar\alpha|\bar\beta)
 \]
for all possible~$\alpha,\beta \in Q^{1,3}$
yields:
\begin{align*}
    r(\{1\},\{1\})&=(-1)^2 T(\{1\}|\{1\})T(\{2,3\}|\{2,3\})=0,\\
    r(\{1\},\{2\})&=(-1)^3 T(\{1\}|\{2\})T(\{2,3\}|\{1,3\})=-12,\\
    r(\{1\},\{3\})&=(-1)^4 T(\{1\}|\{3\})T(\{2,3\}|\{1,2\})=-12,\\
    r(\{2\},\{1\})&=(-1)^3  T(\{2\}|\{1\})T(\{1,3\}|\{2,3\})=-18,\\
    r(\{2\},\{2\})&=(-1)^4  T(\{2\}|\{2\})T(\{1,3\}|\{1,3\})=-6,\\
    r(\{2\},\{3\})&=(-1)^5  T(\{2\}|\{3\})T(\{1,3\}|\{1,2\})=0,\\
    r(\{3\},\{1\})&=(-1)^4 T(\{3\}|\{1\})T(\{1,2\}|\{2,3\})=-6,\\
    r(\{3\},\{2\})&=(-1)^5 T(\{3\}|\{2\})T(\{1,2\}|\{1,3\})=-6,\\
    r(\{3\},\{3\})&=(-1)^6 T(\{3\}|\{3\})T(\{1,2\}|\{1,2\})=-12.
\end{align*}
Thus,  condition~\eqref{eq:tadjugate} holds, and  
Theorem~\ref{prop:diag} implies that~$ TD T^{-1} \in \EP$ for any diagonal matrix~$D \in \R^{3 \times 3}$.
\end{Example}

 \begin{IEEEproof}[Proof of Theorem~\ref{prop:diag}]
 Since~$A=T D T^{-1}$,
 $\exp(At)= T \exp(Dt) T^ {-1}$.
 Pick~$k\in\{1,\dots, \lceil \frac{n-1}{2} \rceil \}$ and~$\alpha \in Q^{k,n}$. Then the corresponding principal minor of~$\exp(At)$ satisfies 
 \begin{align*}
 (\exp(At))(\alpha|\alpha) 
= & \sum_{\beta,\gamma\in Q^{k,n}} T(\alpha|\beta)(\exp(Dt))(\beta|\gamma) T^{-1}(\gamma|\alpha)\\
= &
\sum_{\beta \in Q^{k,n}} T(\alpha|\beta)(\exp(Dt))(\beta|\beta) T^{-1}(\beta|\alpha)\\
=  & 
\sum_{\beta \in Q^{k,n}} \left( T(\alpha|\beta)  T^{-1}(\beta|\alpha)\prod_{i\in\beta} \exp(\lambda_i t) \right),
 \end{align*}
 where we used the fact that~$D$ is diagonal. 
 By Jacobi’s identity~\cite[Chapter~0]{matrx_ana},
 \[
 T^{-1}(\beta|\alpha)= s(\alpha,\beta)  T(\bar\alpha|\bar\beta)/ \det(T).
 \]
Note that $s(\alpha,\beta) = s(\bar \alpha,\bar \beta)$.
 Thus,
  \begin{align*}
&(\exp(At)) (\alpha|\alpha)
 \\&= \frac{1}{\det(T)} 
\sum_{\beta \in Q^{k,n}} \Big (  s(\alpha,\beta) T(\alpha|\beta)    
  T(\bar\alpha|\bar\beta) 
\prod_{i\in\beta} \exp(\lambda_i t) \Big ),
 \end{align*}
and~\eqref{eq:tadjugate} implies that~$(\exp(At))(\alpha|\alpha)>0$.
 \end{IEEEproof}

 \begin{Theorem} \label{thm:genjacobi}
 Suppose that~$A\in\R^{n\times n}$ 
 satisfies the following property: for any~$t\geq 0$ there exists an integer~$k>0$ such that
 \be \label{eq:atellcond}
   I+\frac{At}{\ell} \text{ is TN for any } \ell\geq k. 
 \ee
 Then~$A\in\EP$.
 \end{Theorem}
 
 For example, we say that~$A\in\R^{n\times n}$ is a \emph{weak Jacobi matrix} if it is tri-diagonal and all the entries on the super- and sub-diagonals are non-negative. 
If~$A $ is a weak Jacobi matrix then
it satisfies   the property in Theorem~\ref{thm:genjacobi} (see~\cite[p.~6]{total_book}), so~$A\in\EP$.  In particular,  weak Jacobi   matrices are the generators (in the Lie-algebraic sense) of the group of non-singular TN matrices~\cite{Loewner1955}. 
 
 \begin{IEEEproof}
Fix~$t\geq0$.  Recall that
\be\label{eq:expat}
\exp(At)=\lim_{\ell\to\infty}(I+\frac{At}{\ell})^\ell 
\ee
(see e.g.~\cite{hall}). For any~$\ell$ sufficiently large, $I+\frac{At}{\ell}$ is~TN, and using the fact that the product of two TN matrices is a TN matrix implies that~$(I+\frac{At}{\ell})^j$  is TN for any~$j\geq 1$.
 Let~$B:=(I+\frac{At}{\ell})^\ell$. Using the Cauchy-Binet formula gives that any principal minor of~$B$ satisfies 
 \begin{align*}
     B(\alpha|\alpha)&\geq \left (  (I+\frac{At}{\ell})(\alpha|\alpha) \right )^\ell,
 \end{align*}
 so
 \begin{align*}
     B(\alpha|\alpha)&\geq \left (  1+ \frac{t}{\ell}\sum_{i\in \alpha} a_{ii}+o(\frac{1}{\ell}) \right )^\ell ,
 \end{align*}
 and taking~$\ell \to \infty$ gives 
  \begin{align*}
     (\exp  (At) )(\alpha|\alpha) 
     &\geq \exp( t  \sum_{i\in \alpha} a_{ii})>0,
  \end{align*}
so $A  \in \EP$.
 \end{IEEEproof}

 \begin{Theorem}  \label{thm:signsym}
Let $A \in \R^{n \times n}$. If there exists an~$\bar \varepsilon>0$ such that for any~$ \varepsilon \in[0,\bar \varepsilon]  $ and any integer~$\ell>0$,
\be \label{eq:signsym}
(I + \varepsilon A)^\ell   \text{ is sign symmetric} .
\ee
Then $A \in \EP$. 
 \end{Theorem}
 \begin{IEEEproof}
 Suppose that~$B\in\R^{n\times n}$ is a P-matrix and that~$B^\ell$ is
sign symmetric for any integer~$\ell>0$. 
 Pick $k \in \{1, \dots, n\}$ and~$\alpha \in Q^{k,n}$.
 Then
the Cauchy-Binet formula gives:
\begin{align*}
    (B^2)(\alpha | \alpha ) &= \sum_{\gamma \in Q^{k,n}} B(\alpha | \gamma) B(\gamma | \alpha ) \\
    &\geq (B(\alpha | \alpha))^2\\&>0 ,
\end{align*}
so~$B^2$ is also a  P-matrix and, by assumption, also sign symmetric. Continuing in this fashion implies that
\be\label{eq:citf}
(B^{2^\ell})(\alpha|\alpha) \geq (B(\alpha|\alpha))^{2^\ell } \text{ for any integer } \ell>0.
\ee

Now suppose that~\eqref{eq:signsym}  holds.
Fix~$t\geq 0$. 
Then  for any~$\ell>0$ sufficiently large the matrix~$B:= I+\frac{At}{2^\ell} $ 
is a P-matrix and~\eqref{eq:signsym}
implies that~$B^\ell$ is sign-symmetric, so~\eqref{eq:citf} gives
\begin{align*}
    ( (I + \frac{At}{2^\ell})^{2^\ell}) (\alpha|\alpha) & \geq (  (I+ \frac{At}{2^\ell})  (\alpha|\alpha)  )^{2^\ell}\\
    &\geq (  1+\frac{t}{2^\ell} \sum_{i\in\alpha} a_{ii}+o(\frac{1}{2^\ell} )   )^{2^\ell}.
\end{align*}
Taking~$\ell\to\infty$ gives 
\[
(\exp(At))(\alpha|\alpha)\geq \exp(t\sum_{i\in \alpha} a_{ii} )>0,
\]
which completes the proof. 
 \end{IEEEproof}

Using the relation between sign-symmetric matrices and sign-pattern symmetric matrices, we can restate  Theorem~\ref{thm:signsym} as follows. 

\begin{Corollary} \label{coro:newver}
Suppose that $A \in \R^{n \times n}$ satisfies the following property: there exists an~$\bar \varepsilon>0$ such that 
$ ((I + \varepsilon A)^{(k)}) ^\ell $  is sign-pattern  symmetric for any $  \varepsilon \in[0,\bar \varepsilon],\; 
k\in\{1,\dots, n-1\}$,  and any integer $\ell>0$. 
Then $A \in \EP$.
\end{Corollary}

\begin{Remark}
It follows from~\eqref{eq:poyrt} that 
\[
((I + \varepsilon A)^{(k)}) ^\ell = I+\varepsilon \ell A^{[k]} +o(\varepsilon), 
\]
so the condition in the corollary
   implies that 
 \be\label{eq:akss}
 A^{[k]} \text{ is sign-pattern symmetric for any } k \in \{1, \dots, n-1\}.
 \ee
\end{Remark}


The next example demonstrates that Theorem~\ref{thm:signsym} is indeed more general than Theorem~\ref{thm:epmc}.

\begin{Example}
Consider again the matrix~$A$ in~\eqref{eq:adoesnt}.
Note that $A$ does not belong to any of the  three matrix
classes   in Theorem~\ref{thm:epmc}. A computation yields
\begin{align*}
(I + \varepsilon A)^\ell =& \frac{1}{3}
\begin{bmatrix}
* &  a_1^\ell - a_2^\ell    & * \\
 2 a_1^\ell - 2 a_2^\ell  & * & * \\
0 & 0 & * 
\end{bmatrix}, \\
\left( (I + \varepsilon A)^{(2)} \right)^\ell =& \frac{1}{3}
\begin{bmatrix}
* & * & * \\
0 & * & 
 a_3^\ell - a_4^\ell  \\
0 &  2 a_3^\ell - 2 a_4^\ell  & * 
\end{bmatrix},
\end{align*}
where $a_1 := 1-3 \varepsilon $, $a_2:=  1 - 6 \varepsilon$, $a_3:= 1 -12 \varepsilon + 27 \varepsilon^2$, and $a_4:= 1 -15 \varepsilon + 54 \varepsilon^2$.
Note that both matrices are sign-pattern symmetric. By Corollary~\ref{coro:newver}, $A\in\EP$.
\end{Example}

The next example shows that~\eqref{eq:signsym} is not a necessary condition for~$\EP$.

\begin{Example} \label{exa:Tdiag1}
Consider the matrix
\be \label{eq:Tdiag1}
A = 
\begin{bmatrix}
15 & -9 & -18 \\
1 & 3 & -12 \\
-1 & 3 & 12
\end{bmatrix}.
\ee
Let $D = \diag(0, 12, 18)$ and consider the matrix $T \in \R^{3 \times 3}$   given in Example~\ref{exa:Tdiag}. Then $A =  T D T^{-1}$. Hence, Theorem~\ref{prop:diag} implies that $A \in \EP$. However, $A$ is not sign-pattern symmetric, that is,~condition~\eqref{eq:signsym} is not satisfied.
\end{Example}

\begin{Remark}
Recall that any principal submatrix of a P-matrix is also a P-matrix. A similar conclusion generally does not hold for the $\EP$ case. For example, the matrix $A$ in~\eqref{eq:Tdiag1} is~$\EP$ and has the principal submatrix $A[\{1,2\} | \{1,2\} ] = \begin{bmatrix} 15 & -9 \\ 1 & 3 \end{bmatrix}$, which is not $\EP$ according to Proposition~\ref{prop:case2by2}. 
\end{Remark}

\section{An Application to  Opinion Dynamics}\label{sec:app}
Opinion dynamics  studies how 
local social interactions between social actors lead to the formation of   opinions. Applications include
group decision making, the propagation of rumors, successful
marketing, the emergence of extremism,
and reaching  consensus. For a recent survey, see e.g.,~\cite{survey_opi_dyn}.

Opinion dynamics models   can be classified 
into several  groups. There are
models where the opinions are considered discrete (often accepting only two different
results, say, voting for  
the Republicans or the  Democrats). Examples include
using the Ising model to
describe the behavior of laborers in a strike, where the two options are either to work or to strike~\cite{galam1982}, 
and other models for binary-state dynamics on a network like the Sznajd model~\cite{sznajd2000}.

A second group of models considers  opinions   as 
 continuous variables   that can take values in an interval. The value of the variables may represent 
 the worthiness of a choice.
Examples include the well-known 
Deffuant model~\cite{Deffuant2000}
and the Hegselmann-Krause~\cite{Hegsekmann_2002}  models.  

A third group of models 
considers opinions that are observed as discrete actions, but are
represented internally by each agent as a continuous variable. A typical example is to model the opinion of agent~$i$ as a variable~$x_i$ taking values in~$[-1,1]$,   such that~$x_i>0$ is interpreted as voting for one candidate and~$x_i<0$ represents the voting for the other candidate~\cite{cont_opinion_dyn}.

 Many algorithms have been presented   
  to solve distributed problems with many cooperating agents, such   as average consensus,
rendezvous, and sensor coverage~\cite{1470239,coverge,Rendezvous,eger10}.

Here, we consider the well-known  linear 
consensus algorithm
\be\label{eq:cons}
\dot x=-L x,
\ee
where~$L\in\R^{n\times n}$ is a Laplacian matrix. 
 Note that these  algorithms have also appeared in  many other domains, e.g., 
 average consensus,
rendezvous, and sensor coverage~\cite{1470239,coverge,Rendezvous,eger10}.

We introduce the following requirement. 

\begin{Definition}\label{def:nons}
We say that~\eqref{eq:cons} reaches consensus \emph{without sign-reversal} if for any~$x(0)\in\R^n \setminus \{0\}$ the following properties hold:
\begin{enumerate}
    \item The solution~$x(t)$ converges,  as~$t\to \infty$, to $c 1_n$ for some~$c\ \in \R$;
    \item \label{item:nsr} for any~$t\geq 0$ there exists an index~$i\in\{1,\dots,n\}$ (that may depend on~$x(0)$ and~$t$) such that~$x_i(t) x_i(0) >  0$. 
\end{enumerate}
\end{Definition}
The second requirement aims to prevent a situation where there is a sign-reversal between the vectors~$x(0)$ and~$x(t)$. For example, if~$n=3$ and~$x(0)$ has the sign pattern~$\begin{bmatrix} +& - &+ \end{bmatrix}^T$ we do not allow~$x(t)=\begin{bmatrix}
-&+&-
\end{bmatrix}^T$ at any time~$t\geq 0$. This has a natural interpretation. We do not  allow that \emph{all} the agents     ``change their mind'' along the course of reaching consensus. 
This new requirement is related to a well-known notion in opinion dynamics called ``stubborn agents''~\cite{stubborn}, that is,  an
agent that does not change its opinion
over time (the terminology in this field is not uniform and such agents are also called a leader~\cite{leader2019},
social media~\cite{media2020}, closed-minded~\cite{closed-minded2017},
and  inflexible agent~\cite{inflexible2019}). 
However, the requirement in Definition~\ref{def:nons} is novel, as it does not pinpoint a specific agent, but rather focuses on the global behaviour.

The conditions needed for
reaching consensus are well-known, so we focus on the second requirement in Definition~\ref{def:nons}. 
The next result follows immediately from the fact that the solution of~\eqref{eq:cons} is~$x(t)=\exp(-L t) x(0)$. 
\begin{Proposition}
The second requirement in Definition~\ref{def:nons} holds iff~$ L\in \EP$.   
\end{Proposition}

If~$L$ represents an undirected connectivity  graph then it is symmetric, and Theorem~\ref{thm:epmc} implies that~$ L\in \EP$. But a non-symmetric Laplacian matrix~$L$ may   have complex  eigenvalues and then Theorem~\ref{thm:complex_eig} implies that~$ L\not \in \EP$, so there exists~$t>0$ such that~$\exp(-Lt) \not \in \P$, and by   Theorem~\ref{thm:pmat} 
there exists~$x(0)\in\R^n \setminus\{0\}$
such that~$x_i(0)x_i(t)\leq 0$ for all~$i$. 
 
\subsection{Discrete-time consensus}
Consider the discrete-time linear consensus system 
\[
x(j+1) =  A x(j),
\]
and assume that we require that for any~$x(0)\in\R^n\setminus\{0\}$  the following property holds:
for any time~$k=0,1,\dots$ there exists an index~$i$ (that may depend on~$x(0)$ and~$k$)
such that 
\be\label{eq:reypu}
x_i(k)x_i(0)>0.
\ee
This is the discrete-time analogue 
of condition~\ref{item:nsr}) in Definition~\ref{def:nons}. Since~$x(k)=A^k x(0)$,
\eqref{eq:reypu} is equivalent to 
the requirement that~$A^k\in\P$.
Recall that a matrix $A \in \R^{n \times n}$ is called a PM-matrix if~$A^k \in \P$ for all~$k=0,1,\dots$ (see, e.g., \cite{hershkowitz1986spectra}).  In this sense, EP-matrices can be viewed as the  continuous-time analogue of PM-matrices.

\section{Conclusion} 
P-matrices play an important role in many fields of applied mathematics. Here, we defined the new  notion of an EP-matrix i.e., a matrix~$A$ such that~$\exp(At)$ is a P-matrix for all~$t\geq 0$ (and then also for all~$t\leq 0$).
We show that EP-matrices must  have all real eigenvalues and provide several conditions guaranteeing that a matrix is~EP. 
For $A \in \R^{2\times 2}$, a simple sign pattern condition on the matrix entries is proved to be necessary and sufficient  for~$A \in \EP$.  Using the sign non-reversal property of P-matrices, we also described a natural  application of~$\EP$ matrices to linear  consensus algorithms. 

An interesting topic  for further research is to identify sign patterns that guarantee that a matrix is~EP. Several authors studied dynamical systems that include complementarity  constraints~\cite{OPT_CONT_COMP_SYST,Complementarity_hybrid}, and  the  notion of EP-matrices may also find applications in such systems.

\section*{Acknowledgements}
The authors are grateful to J\"{u}ergen Garloff and Daniel Zelazo for helpful discussions. 

\bibliographystyle{IEEEtranS}
\bibliography{pmatrix}

\begin{thebibliography}{10}
\providecommand{\url}[1]{#1}
\csname url@samestyle\endcsname
\providecommand{\newblock}{\relax}
\providecommand{\bibinfo}[2]{#2}
\providecommand{\BIBentrySTDinterwordspacing}{\spaceskip=0pt\relax}
\providecommand{\BIBentryALTinterwordstretchfactor}{4}
\providecommand{\BIBentryALTinterwordspacing}{\spaceskip=\fontdimen2\font plus
\BIBentryALTinterwordstretchfactor\fontdimen3\font minus
  \fontdimen4\font\relax}
\providecommand{\BIBforeignlanguage}[2]{{%
\expandafter\ifx\csname l@#1\endcsname\relax
\typeout{** WARNING: IEEEtranS.bst: No hyphenation pattern has been}%
\typeout{** loaded for the language `#1'. Using the pattern for}%
\typeout{** the default language instead.}%
\else
\language=\csname l@#1\endcsname
\fi
#2}}
\providecommand{\BIBdecl}{\relax}
\BIBdecl

\bibitem{box_inv2009}
A.~Abate, A.~Tiwari, and S.~Sastry, ``Box invariance in biologically-inspired
  dynamical systems,'' \emph{Automatica}, vol.~45, no.~7, pp. 1601--1610, 2009.

\bibitem{banaji2007}
M.~Banaji, P.~Donnell, and S.~Baigent, ``P matrix properties, injectivity, and
  stability in chemical reaction systems,'' \emph{SIAM J. Applied Math.},
  vol.~67, no.~6, pp. 1523--1547, 2007.

\bibitem{comp_barshalom_omri}
\BIBentryALTinterwordspacing
E.~Bar-Shalom, O.~Dalin, and M.~Margaliot, ``Compound matrices in systems and
  control theory: a tutorial,'' 2022, submitted. [Online]. Available:
  \url{https://arxiv.org/abs/2204.00676}
\BIBentrySTDinterwordspacing

\bibitem{media2020}
H.~Z. Brooks and M.~A. Porter, ``A model for the influence of media on the
  ideology of content in online social networks,'' \emph{Phys. Rev. Research},
  vol.~2, p. 023041, 2020.

\bibitem{closed-minded2017}
B.~Chazelle and C.~Wang, ``Inertial {Hegselmann–Krause} systems,'' \emph{IEEE
  Trans.\ Automat.\ Control}, vol.~62, no.~8, pp. 3905--3913, 2017.

\bibitem{choudhury2021}
\BIBentryALTinterwordspacing
P.~N. Choudhury, ``Characterizing total positivity: single vector tests via
  linear complementarity, sign non-reversal, and variation diminution,'' 2021.
  [Online]. Available: \url{https://arxiv.org/abs/2103.05624}
\BIBentrySTDinterwordspacing

\bibitem{Apoorva2021}
P.~N. Choudhury, M.~R. Kannan, and A.~Khare, ``Sign non-reversal property for
  totally non-negative and totally positive matrices, and testing total
  positivity of their interval hull,'' \emph{Bull. London Math. Soc.}, pp.
  1--10, 2021.

\bibitem{coverge}
J.~Cortes, S.~Martinez, T.~Karatas, and F.~Bullo, ``Coverage control for mobile
  sensing networks,'' \emph{IEEE Trans. Robotics and Automation}, vol.~20,
  no.~2, pp. 243--255, 2004.

\bibitem{LCP_cottle}
R.~W. Cottle, J.-S. Pang, and R.~E. Stone, \emph{The Linear Complementarity
  Problem}.\hskip 1em plus 0.5em minus 0.4em\relax Society for Industrial and
  Applied Mathematics, 2009.

\bibitem{coxson}
G.~E. Coxson, ``The {P-matrix} problem is {co-NP-complete},''
  \emph{Mathematical Programming}, vol.~64, pp. 173--178, 1994.

\bibitem{Deffuant2000}
G.~Deffuant, D.~Neau, F.~Amblard, and G.~Weisbuch, ``Mixing beliefs among
  interacting agents,'' \emph{Advances in Complex Systems}, vol.~03, no. 01n04,
  pp. 87--98, 2000.

\bibitem{total_book}
S.~M. Fallat and C.~R. Johnson, \emph{Totally Nonnegative Matrices}.\hskip 1em
  plus 0.5em minus 0.4em\relax Princeton, NJ: Princeton University Press, 2011.

\bibitem{fallat2017total}
S.~Fallat, C.~R. Johnson, and A.~D. Sokal, ``Total positivity of sums,
  {H}adamard products and {H}adamard powers: Results and counterexamples,''
  \emph{Linear Algebra Appl.}, vol. 520, pp. 242--259, 2017.

\bibitem{FIEDLER1966}
M.~Fiedler and V.~Ptak, ``Some generalizations of positive definiteness and
  monotonicity,'' \emph{Numerische Mathematik}, vol.~9, pp. 163--172, 1966.

\bibitem{fiedler1962matrices}
------, ``On matrices with non-positive off-diagonal elements and positive
  principal minors,'' \emph{Czechoslovak Mathematical Journal}, vol.~12, no.~3,
  pp. 382--400, 1962.

\bibitem{galam1982}
S.~Galam, Y.~{Gefen (Feigenblat)}, and Y.~Shapir, ``Sociophysics: A new
  approach of sociological collective behaviour: {I.} mean-behaviour
  description of a strike,'' \emph{J. Math. Sociology}, vol.~9, no.~1, pp.
  1--13, 1982.

\bibitem{gale-nikaido1965}
D.~Gale and H.~Nikaido, ``The {Jacobian} matrix and global univalence of
  mapping,'' \emph{Mathematische Annalen}, vol. 159, p. 81–93, 1965.

\bibitem{gk_book}
F.~R. Gantmacher and M.~G. Krein, \emph{Oscillation Matrices and Kernels and
  Small Vibrations of Mechanical Systems}.\hskip 1em plus 0.5em minus
  0.4em\relax Providence, RI: American Mathematical Society, 2002, translation
  based on the~1941 {Russian} original.

\bibitem{gen_Gale_2}
C.~B. Garcia and W.~I. Zangwill, ``On univalence and {P-matrices},''
  \emph{Linear Algebra Appl.}, vol.~24, pp. 239--250, 1979.

\bibitem{hall}
B.~C. Hall, \emph{Lie Groups, Lie Algebras, and Representations: An Elementary
  Introduction}.\hskip 1em plus 0.5em minus 0.4em\relax Springer, 2003.

\bibitem{Hegsekmann_2002}
R.~Hegselmann and U.~Krause, ``Opinion dynamics and bounded confidence: models,
  analysis and simulation,'' \emph{J. Artificial Societies and Social
  Simulation}, vol.~5, no.~3, 2002.

\bibitem{hershkowitz1986spectra}
D.~Hershkowitz and C.~R. Johnson, ``Spectra of matrices with {P}-matrix
  powers,'' \emph{Linear Algebra Appl.}, vol.~80, pp. 159--171, 1986.

\bibitem{hersh2003}
D.~Hershkowitz and N.~Keller, ``Positivity of principal minors, sign symmetry
  and stability,'' \emph{Linear Algebra Appl.}, vol. 364, pp. 105--124, 2003.

\bibitem{sigmund_evolution_book}
J.~Hofbauer and K.~Sigmund:, \emph{The Theory of Evolution and Dynamical
  Systems}.\hskip 1em plus 0.5em minus 0.4em\relax Cambridge University Press,
  1988.

\bibitem{topics_math_ana}
R.~A. Horn and C.~R. Johnson, \emph{Topics in Matrix Analysis}.\hskip 1em plus
  0.5em minus 0.4em\relax Cambridge University Press, 1991.

\bibitem{matrx_ana}
------, \emph{Matrix Analysis}, 2nd~ed.\hskip 1em plus 0.5em minus 0.4em\relax
  Cambridge University Press, 2013.

\bibitem{inflexible2019}
F.~Jacobs and S.~Galam, ``Two-opinions-dynamics generated by inflexibles and
  non-contrarian and contrarian floaters,'' \emph{Adv. Complex Syst.}, vol.~22,
  no.~04, p. 1950008, 2019.

\bibitem{mat_posi_johnson}
C.~R. Johnson, R.~L. Smith, and M.~J. Tsatsomeros, \emph{Matrix Positivity},
  ser. Cambridge Tracts in Mathematics.\hskip 1em plus 0.5em minus 0.4em\relax
  Cambridge University Press, 2020, vol. 221.

\bibitem{diag_stab_book}
E.~Kaszkurewicz and A.~Bhaya, \emph{Matrix Diagonal Stability in Systems and
  Computation}.\hskip 1em plus 0.5em minus 0.4em\relax Basel: Birkhauser, 2000.

\bibitem{rami_osci}
R.~Katz, M.~Margaliot, and E.~Fridman, ``Entrainment to subharmonic
  trajectories in oscillatory discrete-time systems,'' \emph{Automatica}, vol.
  116, p. 108919, 2020.

\bibitem{kellogg1972complex}
R.~Kellogg, ``On complex eigenvalues of {M}- and {P}-matrices,''
  \emph{Numerische Mathematik}, vol.~19, no.~2, pp. 170--175, 1972.

\bibitem{Rendezvous}
J.~Lin, A.~S. Morse, and B.~D.~O. Anderson, ``The multi-agent rendezvous
  problem. {Part} 1: The synchronous case,'' \emph{SIAM J.\ Control Optim.},
  vol.~46, no.~6, pp. 2096--2119, 2007.

\bibitem{Loewner1955}
C.~Loewner, ``On totally positive matrices,'' \emph{Mathematische Zeitschrift},
  vol.~63, no.~1, pp. 338--340, 1955.

\bibitem{margaliot2019revisiting}
M.~Margaliot and E.~D. Sontag, ``Revisiting totally positive differential
  systems: A tutorial and new results,'' \emph{Automatica}, vol. 101, pp.
  1--14, 2019.

\bibitem{cont_opinion_dyn}
A.~C.~R. Martins, ``Continuous opinions and discrete actions in opinion
  dynamics problems,'' \emph{Int. J. Modern Physics C}, vol.~19, no.~04, pp.
  617--624, 2008.

\bibitem{gen_Gale_1}
A.~Mas-Colell, ``Homeomorphisms of compact, convex sets and the {Jacobian}
  matrix,'' \emph{SIAM J. Math. Anal.}, vol.~10, no.~6, pp. 1105--1109, 1979.

\bibitem{eger10}
M.~Mesbahi and M.~Egerstedt, \emph{Graph-Theoretic Methods in Multiagent
  Networks}.\hskip 1em plus 0.5em minus 0.4em\relax Princeton, NJ: Princeton
  Univ. Press, 2010.

\bibitem{MOYLAN197753}
P.~J. Moylan, ``Matrices with positive principal minors,'' \emph{Linear Algebra
  Appl.}, vol.~17, no.~1, pp. 53--58, 1977.

\bibitem{muldo1990}
J.~S. Muldowney, ``Compound matrices and ordinary differential equations,''
  \emph{The Rocky Mountain J. Math.}, vol.~20, no.~4, pp. 857--872, 1990.

\bibitem{murty_pmat}
K.~G. Murty, ``On the number of solutions to the complementarity problem and
  spanning properties of complementary cones,'' \emph{Linear Algebra Appl.},
  vol.~5, p. 65–108, 1972.

\bibitem{murty_LCP}
------, \emph{Linear Complementarity, Linear and Nonlinear Programming}.\hskip
  1em plus 0.5em minus 0.4em\relax Berlin: Heldermann, 1988.

\bibitem{nika_book}
H.~Nikaido, \emph{Convex Structures and Economic Theory}.\hskip 1em plus 0.5em
  minus 0.4em\relax Academic Press, 1968.

\bibitem{survey_opi_dyn}
H.~Noorazar, ``Recent advances in opinion propagation dynamics: a 2020
  survey,'' \emph{Eur. Phys. J. Plus}, vol. 135, p. 521, 2020.

\bibitem{bmatrices}
J.~M. Pena, ``A class of {P-matrices} with applications to the localization of
  the eigenvalues of a real matrix,'' \emph{SIAM J. Matrix Anal. Appl.},
  vol.~22, pp. 1027--1037, 2001.

\bibitem{pinkus}
A.~Pinkus, \emph{Totally Positive Matrices}.\hskip 1em plus 0.5em minus
  0.4em\relax Cambridge, UK: Cambridge University Press, 2010.

\bibitem{posi-tutorial}
A.~Rantzer and M.~E. Valcher, ``A tutorial on positive systems and large scale
  control,'' in \emph{{Proc.\ 57th IEEE Conf. on Decision and Control}}, Miami
  Beach, FL, USA, 2018, pp. 3686--3697.

\bibitem{1470239}
W.~Ren, R.~W. Beard, and E.~M. Atkins, ``A survey of consensus problems in
  multi-agent coordination,'' in \emph{Proc. American Control Conf.}, vol.~3,
  2005, pp. 1859--1864.

\bibitem{schwarz1970}
B.~Schwarz, ``Totally positive differential systems,'' \emph{Pacific J. Math.},
  vol.~32, no.~1, pp. 203--229, 1970.

\bibitem{hlsmith}
H.~L. Smith, \emph{Monotone Dynamical Systems: An Introduction to the Theory of
  Competitive and Cooperative Systems}, ser. Mathematical Surveys and
  Monographs.\hskip 1em plus 0.5em minus 0.4em\relax Providence, RI: Amer.
  Math. Soc., 1995, vol.~41.

\bibitem{sznajd2000}
K.~Sznajd-Weron and J.~Sznjad, ``Opinion evolution in closed community,''
  vol.~11, no.~6, pp. 1157--1165, 2000.

\bibitem{vol_comp}
Y.~Takeuchi and N.~Adachi, ``The existence of globally stable equilibria of
  ecosystems of the generalized {Volterra} type,'' \emph{J. Math. Biology},
  vol.~10, p. 401–415, 1980.

\bibitem{gen_p_matrices}
M.~J. Tsatsomeros, ``Generating and detecting matrices with positive principal
  minors,'' \emph{Asian Information-Science-Life: An International Journal},
  vol.~1, no.~2, pp. 115--132, 2002.

\bibitem{Complementarity_hybrid}
A.~van~der Schaft and J.~Schumacher, ``Complementarity modeling of hybrid
  systems,'' \emph{IEEE Trans.\ Automat.\ Control}, vol.~43, no.~4, pp.
  483--490, 1998.

\bibitem{OPT_CONT_COMP_SYST}
A.~Vieira, B.~Brogliato, and C.~Prieur, ``Quadratic optimal control of linear
  complementarity systems: First-order necessary conditions and numerical
  analysis,'' \emph{IEEE Trans.\ Automat.\ Control}, vol.~65, no.~6, pp.
  2743--2750, 2020.

\bibitem{Hausdorff_contract}
\BIBentryALTinterwordspacing
C.~Wu, R.~Pines, M.~Margaliot, and J.~Slotine, ``Generalization of the
  multiplicative and additive compounds of square matrices and contraction
  theory in the {Hausdorff} dimension,'' \emph{IEEE Trans.\ Automat.\ Control},
  2021, to appear. [Online]. Available: \url{https://arxiv.org/abs/2008.10321}
\BIBentrySTDinterwordspacing

\bibitem{leader2019}
\BIBentryALTinterwordspacing
Y.~Yi and S.~Patterson, ``Disagreement and polarization in two-party social
  networks,'' 2019. [Online]. Available: \url{https://arxiv.org/abs/1911.11338}
\BIBentrySTDinterwordspacing

\bibitem{stubborn}
E.~Yildiz, A.~Ozdaglar, D.~Acemoglu, A.~Saberi, and A.~Scaglione, ``Binary
  opinion dynamics with stubborn agents,'' \emph{ACM Trans. Econ. Comput.},
  vol.~1, no.~4, 2013.

\end{thebibliography}

 \end{document}